\documentclass{article}
\usepackage{epsfig}
\usepackage{amsmath, amssymb}
\usepackage{amsbsy}
\usepackage{graphicx}
\usepackage{changebar}
\numberwithin{equation}{section}
\begin{document}

\centerline{FRACALMO PRE-PRINT   {\bf www.fracalmo.org}}
\centerline{Paper  published in}
\centerline{\bf   Vietnam Journal of Mathematics,  Vol. 32 (SI),  65-75 (2004)}
\vskip 0.5truecm
%% \cen{Springer Verlag:  ISSN 0866-7179} \vsh
\hrule

\vskip 1truecm

\font\title=cmbx12 scaled\magstep1
\font\bfs=cmbx12 scaled\magstep1

%% \vsh
\begin{center}
{\title{From Power Laws to Fractional Diffusion:}}

\medskip

{\title{the Direct Way}}

\vskip 0.5truecm

{Rudolf GORENFLO}$^{(1)}$
 and {Entsar A.A ~ABDEL-REHIM} $^{(2)}$

%%%%%%%%%%%%%%%%

\vskip 0.25truecm

$\null^{(1)}$
 Fachbereich Mathematik \& Informatik \\
Erstes Mathematisches Institut,
 Freie   Universit\"at  Berlin, \\
 Arnimallee  3, D-14195 Berlin, Germany \\
%% Tel: +49-30-838.75.352 $\;$ Fax: +49-30-838.75.403 $\;$ \\
E-mail: {\tt gorenflo@mi.fu-berlin.de}
\\ [0.25 truecm]

$\null^{(2)}$
 On leave from Department of Mathematics and Computer Science, \\
  Faculty of Science, Suez Canal University, Egypt 
% E-mail: {\tt mainardi@bo.infn.it}
\\ [0.25 truecm]

\end{center}
\begin{abstract}
Starting from the model of continuous time random walk (Montroll and
Weiss 1965) that can also be considered as a compound renewal process
we focus our interest on random walks in which the probability
distributions of the waiting times and jumps  have fat tails
characterized by power laws with exponent between 0 and 1 for the
waiting times, between 0 and 2  for the jumps. By stating the
relevant lemmata (of Tauber type) for the distribution functions  we need not
distinguish between continuous and  discrete space and time. We will
see that by a well-scaled passage to the 
diffusion limit diffusion processes fractional in time as well as in
space are obtained. The corresponding equation of evolution is a
linear partial pseudo-differential equation with fractional
derivatives in time and in space, the orders being equal to the above
exponents. Such processes are enjoying increasing popularity in
applications in physics, chemistry, finance and other fields, and
their behaviour can be well approximated and visualized by simulation
via various types of  random walks. For their explicit solutions there
are available integral representations that allow to investigate their
detailed structure. For ease of presentation we restrict attention to
the spatially one-dimensional symmetric situation.  
\\
\\
\textbf{MSC 2000}: 26A33, 33E12, 45E10, 45K05, 60F05, 60G50, 60J60 
\\
\textbf{Keywords}: continuous time random walks, asymptotics of
Fourier and Laplace transforms, convergence in law, well-scaled
passage to the diffusion limit, Mittag-Leffler function.

\end{abstract}
%%%%%%%%%%%%%%%%%%%%%%%%%%%%%%%%%%%%%%%%
\newpage
%%%%%%%%%%%
\section{ Introduction: concepts and notations }

    We consider spatially one-dimensional random walks of the
    following basic structure. A particle (or a wanderer) starting at
   the  time instant $t=0$   at the space point $x=0$  makes jumps of
    random size $X_{k} $
    in random instants $t_{k}$, $k \in \mathbb{N}=\{1,2,3,\cdots\}$,
    $0<t_{1}<t_{2}\cdots \to \infty $. For convenience we set
    $t_{0}=0$. Then in 
    the time interval $t_{n}\le t<t_{n+1}$   the particle is sitting
    in the point $x= S_{n}:=\sum\limits_{k=1} ^{n} X_{k}$. We
    assume the jumps to be \emph{independent identically distributed
    (iid)} random variables, all having the same probability
    distribution as a generic real random variable $X$ , called the
    jump. Likewise we assume the waiting times $T_{k}:=t_{k}-
    t_{k-1}$   to be  iid  random 
    variables, all equal in distribution to a generic non-negative
    random variable $T$ , called the waiting time. So, this process is
    what in mathematical literature is called a \emph{compound} (or
    \emph{cumulative}) renewal process \cite{C}. We denote  the
    distribution 
    functions of the waiting time $T$  and the jump $X$  by
    $\Phi$   and  $W$, respectively, by 
\[P(T\le t)= \Phi(t)\;,\; 0\le t<\infty \; ;\;\; P(X\le x)= W(x)\; ,\;
    -\infty<x<\infty \; .\]
      Conveniently using the language of \emph{generalized
     functions} in the sense of \cite{GS} or \cite{Z} we introduce the
     (generalized) probability \emph{densities}  $\phi$  and $w$, so that    
\[ 
\Phi(t)=  \int\limits_{0}^{t}\phi(t') dt'
  \;,0 \le t < \infty \;; \; W(x)= \int\limits_{-\infty}^{x}w(x')dx'
  \;, -\infty <x<\infty\;.
\]
Denoting for time  instant $t$ the probability density to find the particle in
point $x$    by $p(x,t)$    we then have, by conditioning on the last
jump before  $t$ and using the delta function $\delta(x)$, for $ 0\le
t <\infty  $ and $-\infty<x<\infty$  the  integral equation (see
\cite{ MW}) of \emph{continuous time  random walk}
\begin{equation}\label{E:eq1}
p(x,t)=\delta(x)(1-\Phi(t))\,+\,\int\limits_{0}^{t}
\{
\int\limits_{-\infty}^{\infty} w(x-x')
 p(x',t') dx' \}  \,\phi(t-t') dt' \; .
\end{equation}
 For the cumulative function  $P(x,t)=
\int\limits_{-\infty}^{x} p(x',t)\, dx'$ we have, with the Heaviside
step function $H(x)$,   the equation  
\begin{equation}
P(x,t)= H(x)(1-\Phi(t))\,+\,\int\limits_{0}^{t}
\{
\int\limits_{-\infty}^{\infty} W(x-x')
 dP(x',t') \} \,d\Phi(t-t') \;.
\end{equation}
 To proceed further we use the machinery of the transforms of Laplace
 and Fourier. The general formulas for $s\ge 0$, $-\infty <\kappa
 <\infty $  are  
\[ \widetilde{g}(s)= \int\limits_{0}^{\infty} e^{-st} g(t) dt=
 \int\limits_{0}^{\infty}e^{-st} d G(t)  \; ; \; \widehat{f}(\kappa)=\int\limits_{-\infty}^{\infty}
e^{i \kappa x} f(x) dx= \int\limits_{-\infty}^{\infty}
e^{i \kappa x} dF(x)  \; .   \]
Essentially applying these formulas to probability densities with $
 0\le t<\infty$ and $-\infty<x<\infty $
 we can safely take the Laplace variable  $s$  as real
and non-negative. Furthermore we will work with convolutions of
(generalized) functions, namely with the Laplace convolution and the
Fourier convolution:  
\[(g_{1}*g_{2})(t)= \int\limits_{0}^{\infty}
g_{1}(t')g_{2}(t-t') dt' \; ; \;
(f_{1}*f_{2})(x)= \int\limits_{-\infty}^{\infty} 
f_{1}(x')f_{2}(x-x') dx'\; . \]
Then, applying the transforms of Fourier and Laplace in succession to
the equation (\ref{E:eq1}) and using the well-known
operational rules, we 
arrive at the relation  
\begin{equation}\label{E:eq3}
\widehat{\widetilde{p}}(\kappa,s)=\frac{1-\widetilde{\phi}(s)}{s}\, + 
\widetilde{\phi}(s)\widehat{w}(\kappa) \widehat{\widetilde{p}}(\kappa,s)  \; ,
\end{equation}
which   leads to the famous Montroll-Weiss equation, see \cite{MW},
\begin{equation}\label{E:eq4}
 \widehat{\widetilde{p}} (\kappa,s) = \frac{1-\tilde{\phi}(s)}{s} \,
\frac{1}{1-\widehat{w}(\kappa) \tilde{\phi}(s)}\;\; .
\end{equation}
This equation can alternatively be derived from the Cox formula, see
    \cite{C} chapter 8 formula (4),  describing the process as
    subordination of a random 
    walk to a renewal process. By inverting the transforms one can, in
    principle, find the evolution $p(x,t)$  of the sojourn density for
    time $t$
    running from zero to infinity.

    Our aim is to show that under appropriate assumptions of power
    laws for the distribution functions $\Phi(t)$, $t \ge 0$,  and
    $W(x)$, $- \infty < x< \infty$, under observance
    of  a scaling relation between  the positive parameters  $h$   and $\tau$
    the re-scaled random walk  $S_{n}(h)= \sum\limits_{k=1}^{n} h
    X_{k} $
 happening at the instants  $t_{n}(\tau)= \sum\limits_{k=1}^{n} \tau
    T_{k}$  (with  $S_{0}(h)=0$, $ t_{0}(\tau)=0 $ ) weakly (or in law)
    tends, for $h$    and  $\tau$  tending
    to zero, to a process obeying the space-time fractional diffusion
    equation. Specifically, we will show that the sojourn probability
    density $p_{h,\tau}(x,t)$   tends weakly  to the
    solution  $u(x,t)$ of the Cauchy problem for $t>0$ and $x \in \mathbb{R}$
\begin{equation}\label{E:eq5}
\underset{t\;\; *}{D}^{\beta}\,u(x,t)=R^{\alpha}\,u(x,t) 
 \; , \; u(x,0)=\delta(x)\; .
\end{equation}
Here $0<\alpha\le 2,\; 0<\beta\le 1$. The fractional Riesz derivative
$R^{\alpha}$ (in space) is defined as follows: the  Fourier
transform of $R^{\alpha} f(x)$ is
$-|\kappa|^{\alpha}\widehat{f}(\kappa)$ for a sufficiently
well-behaved function $f(x)$. Compare \cite{F1952}, \cite{R},
\cite{SKM}. The  Caputo fractional derivative   (in time) can be
defined through its image in the Laplace transform domain. The Laplace
transform of   $\underset{t\;\; *}{D}^{\beta}\,g(t)$ is
$s^{\beta}\widetilde{g}(s) -s^{\beta-1} g(0)$. We have
$
\underset{t\;\; *}{D}^{\beta}\,g(t)=\frac{dg(t)}{dt} \; \text{for}\;
\beta=1$ 
 but 
\[\underset{t\;\; *}{D}^{\beta}\,g(t)=
\frac{1}{\Gamma(1-\beta)}\{\frac{d}{dt} \int\limits_{0}^{t}
(t-t')^{-\beta} g(t') dt'\,-\, t^{-\beta} g(0)
\}\; \text{for}\; 0<\beta<1 \; ,   \]
compare \cite{GM1997}.  
In the Fourier-Laplace domain the Cauchy problem (\ref{E:eq5}) appears
in the form 
$s^{\beta}\widehat{\widetilde{u}} (\kappa,s)- s^{\beta-1}=
-|\kappa|^{\alpha}\widehat{\widetilde{u}} (\kappa,s)   $
from  which
we  obtain 
\begin{equation}\label{E:eq6}
\widehat{\widetilde{u}} (\kappa,s) =\frac{s^{\beta -1}}{s^{\beta}+
  |\kappa|^{\alpha} }\, ,\; s>0\, ,\; \kappa \in \mathbb{R}\; .
\end{equation}

Let us refer to \cite{MLP} for the analytical theory of representing
the function $u(x,t)$  , namely the fundamental solution of the space-time
fractional diffusion equation, in dependence on the parameters
$\alpha$   and $\beta$.

      To carry out the passage to the diffusion limit we state
        in Section $2$  two \emph{MASTER LEMMATA} and four
        simplifications  relating
       the asymptotic 
       behaviours of the distribution functions $W(x)$ and $\Phi(t)$
        near infinity  to the
       asymptotic behaviour of 
       their Laplace and Fourier transforms near zero. Section $ 3$ is
       devoted to the actual passage to the diffusion limit, in Section
       $4$ some 
       examples are presented, and a few  historical comments are given in
       Section $5$.  
%%%%%%%%%%%%*****************************************************
\section{Six lemmata}
\textbf{Definition:}  As in \cite{{BGT}} we call a positive measurable
function $\nu$, defined on some neighbourhood $[x^{*}, \infty)$ of
infinity 
, \emph{slowly varying} if $\nu(ax)/\nu(x) \to 1$  as $x \to \infty$ for every
$a>0$. Examples: $(log x)^{\gamma}$ with $\gamma \in \mathbb{R}$ and
$exp \,(\frac{log x}{log log x})$. 
\\
\\
\textbf{MASTER LEMMA 1:}
 Assume  $W(x)$ increasing, $W(-\infty)=0$, $W(\infty)=1$,
 symmetry
 $\int\limits_{(-\infty,-x)}dW(x')=\int\limits_{(x,\infty)}dW(x')$
   for $x \ge 0$, let $L$ be a slowly varying function  and
 assume either (a) or (b).
\newline
(a)  $\sigma^{2}:=\int\limits_{-\infty}^{\infty} x^{2} dW(x) <\,\infty$,
labelled as  $\alpha=2$ ,
\\
(b) $\int\limits_{(x,\infty)}d W(x) \sim b \alpha^{-1} x^{-\alpha}
L(x)$ 
  for $x \to \infty$,  $\alpha \in (0,2)$  and   $b>0$. 
\\
  Then, with
\begin{equation}\label{E:eq12}
 \mu=\frac{\sigma^{2}}{2}  \; \text{and} \; L(x)\equiv 1\;  \text{ in
   case (a),}\; 
\mu= \frac{b \pi}{\Gamma(\alpha +1) sin (\alpha \pi /2)} \;  \text{ in case
(b)} \; ,
\end{equation}
we have the asymptotics  $1-\widehat{w}(\kappa ) \sim
\mu|\kappa|^{\alpha}L(|\kappa|^{-1})$ for $\kappa \to 0$.
%\end{lemma}
%
\\
\\
\textbf{Comments}
The proof can be distilled from Chapter 8 of \cite{BGT}. In some sense
this lemma is a partial reformulation (with a constant corrected) of
Gnedenko's theorem on the domain of attraction of stable probability
laws, see \cite{GK}. 
\\
\\
\textbf{MASTER LEMMA 2:}
Assume $\Phi(t)$ increasing, $\Phi(0)=0$, $\Phi(\infty)=1$, let $M$ be
a slowly varying function and assume either (A) or (B). 
\\
(A) $\rho:=\int\limits_{0}^{\infty}\, t d \Phi(t)< \infty $, labelled as
$\beta=1$, 
\\
(B) $\int\limits_{(t,\infty)} d \Phi(t)\sim c \beta^{-1} t^{-\beta}
 M(t)$ for $t\to \infty$, $\beta \in (0,1)$ and  $c>0$.
\\
 Then, with
\begin{equation} \label{E:eq13}
 \lambda=\rho \; \text{and}  \; M(t)\equiv 1 \;  \text{in case (A),}
 \; 
\lambda=\frac{c\Gamma(1-\beta)}{\beta}\;   \text{in case (B)} \; ,
\end{equation}
 we have the asymptotics $ 1- \widetilde{\phi}(s) \sim \lambda s^{\beta}
M(s^{-1})$ for $0<s \to 0$. 
%\end{lemma}  
% 
\\
\\
\textbf{Comments}
This lemma is a special case of Karamata's theorem of $1931$. A proof
can be found in the book \cite{BGT}. In the case of existing
non-generalized functions as densities 
     the Master Lemmata imply the two following (convenient) C-Lemmata (by
     integration it can be shown that the assumptions of the Master
     Lemmata are met). See \cite{GM2002}. For the fully equidistant
     discrete case 
     the two D-Lemmata are useful. See \cite{GV}. 
\\
\\
\textbf{C-Lemma $1$} (for jump densities): 
Assume $w(x)\ge 0$, $w(x)=w(-x)$ for $x\in \mathbb{R}$,
$\int\limits_{-\infty}^{\infty} w(x)dx=1$ and either (a) or (b).  
\\
(a) $\sigma^{2}:=\int\limits_{-\infty}^{\infty} x^{2} w(x) dx
<\,\infty$, labelled as $\alpha=2$,
\\
(b) $w(x)\sim b|x|^{-(\alpha+1)}$ for $|x| \to \infty$, 
$\alpha \in (0,2)$ and $b>0$.
\\
 Then  with  $\mu$ as in (\ref{E:eq12}), we have the asymptotics 
\begin{equation}\label{E:eq14}
1- \widehat{w}(\kappa)\sim \mu |\kappa|^{\alpha} \; \text{for} \;
\kappa \to 0 \; .
\end{equation}
%
%\end{proof}
%
%
\textbf{C-Lemma $2$} (for waiting time densities): 
Assume $\phi(t)\ge 0 $ for $t>0$, $\int\limits_{0}^{\infty}\phi(t)
dt=1 $, and either  (A) or (B).
\\
(A) $\rho:=\int\limits_{0}^{\infty}\, t\phi(t) dt<\infty$,  labelled
as $\beta=1$,
\\
(B) $ \phi(t)\sim \, c t^{-(\beta +1)}\; \text{for} \; t \to \infty$, 
 $\beta \in (0,1)$ and $c>0$.
\\
 Then with $\lambda$ as in (\ref{E:eq13}), we have the asymptotics 
\begin{equation}\label{E:eq15}
1- \widetilde{\phi}(s) \sim \lambda s^{\beta} \; .
\end{equation}
%
% \end{proof}
%
\\
\textbf{D-Lemma $1$: } Assume $p_{k}\ge 0$,
    $\sum\limits_{-\infty}^{\infty} p_{k}=1$, 
 symmetry  $p_{k}=p_{-k}$ for all integers  $k\in \mathbb{Z}$, 
 and either (a) or (b).  
\\
(a) $\sigma^{2}:=\sum\limits_{-\infty}^{\infty} k^{2} p_{k}
<\,\infty$, labelled as $\alpha=2$,
\\
(b) $p_{k}\sim b|k|^{-(\alpha+1)}$ for $|k| \to \infty$,
$\alpha \in (0,2)$ and $b>0$.
\\
Then with
    $\mu$ as in (\ref{E:eq12}) we have the asymptotics 
(\ref{E:eq14}).
%\end{lemma}
\\
\\
\textbf{D-Lemma $2$:} Assume $c_{n}\ge 0$,
    $\sum\limits_{n=1}^{\infty}c_{n}=1 $ and either (A) or (B).  
\\
(A) $\rho:=\sum\limits_{1}^{\infty}\, nc_{n}<\infty$,  labelled
as $\beta=1$,
\\
(B) $c_{n}\sim \, c n^{-(\beta +1)}\; \text{for} \; n \to \infty$,
 $\beta \in (0,1)$ and $c>0$.
\\
 Then with $\lambda$ as in (\ref{E:eq13}) we have the asymptotics
(\ref{E:eq15}).
%\end{lemma}
%
%%%%%%%%%**********************************************
\section{ Well-scaled passage to the diffusion limit}
     As already indicated in the Introduction, we multiply the jumps
  $X_{k}$ by a factor  $h$, the waiting times $T_{k}$  by a factor
  $\tau$. So, we get a transformed random walk
  $S_{n}(h)=\sum\limits_{k=1}^{n} h X_{k}$  with jump instants
  $t_{n}(\tau) =\sum\limits_{k=1}^{n} \tau T_{k}$  that we now
  investigate with the aim of passing to the limit $h\to 0$, $\tau
  \to 0$ under a scaling relation between $h$ and $\tau$  yet to be
  established, assuming that the conditions of Master Lemma $1$ and
  Master Lemma $2$ are fulfilled. As it is convenient to work in the
  Fourier-Laplace domain we note that the density $\phi_{\tau}(t) $
  of the reduced waiting times $\tau T_{k}$   and the density
  $w_{h}(x)$ of the reduced jumps $h X_{k}$  are $\phi_{\tau}(t)=\phi(t/ \tau)/ \tau $,
  $t \ge 0$  ; $w_{h}(x)=w(x/h)/h $, $-\infty <x< \infty$. The
  corresponding transforms are simply $\widetilde{\phi_{\tau}}(s)=
  \widetilde{\phi}(s \tau)$,
  $\widehat{w_{h}}(\kappa)=\widehat{w}(\kappa h)$. We are interested
  in the sojourn probability density $p_{h,\tau}(x,t)$   of the
  particle subject to the transformed random walk. In analogy to the
  Montroll-Weiss equation (\ref{E:eq4}) we get
\begin{equation}\label{E:eq7}
\widehat{\widetilde{p}}_{h,\tau} (\kappa,s) =\frac{1-\widetilde{\phi}_{\tau}(s)}{s}
\frac{1}{1-\widehat{w}_{h}(\kappa) \widetilde{\phi}_{\tau}(s)} \;= 
\frac{1-\widetilde{\phi}(\tau s)}{s}
\frac{1}{1-\widehat{w}(h\kappa) \widetilde{\phi}(s\tau)} \; .
\end{equation}
Considering now $s$    and $\kappa $ \emph{fixed} and $\ne 0$ we find for $h\to 0$,
$\tau \to 0$  from the Master Lemmata (replacing there $\kappa$   by
$\kappa h$,  $s$  by  $s \tau$ )  by a trivial calculation, omitting
asymptotically negligible terms and using the slow variation property
$L(1/ (\kappa h))\sim L(1/h)$, $M(1/(s \tau))\sim
M(1/ \tau)$, the asymptotics (\ref{E:eq8}) with (\ref{E:eq9}).
\begin{equation}\label{E:eq8}
\widehat{\widetilde{p}}_{h,\tau}(\kappa,s)=\frac{\lambda
  \tau^{\beta} s^{\beta-1} M(1/ \tau)}{\mu (h|\kappa|)^{\alpha}L(1/ h)+
\lambda (\tau s)^{\beta}M(1/ \tau)
  } = \frac{s^{\beta
  -1}}{r(h, \tau)|\kappa|^{\alpha} + s^{\beta} } \; ,
\end{equation}
\begin{equation}\label{E:eq9}
 r(h,\tau)= \frac{\mu h^{\alpha}L(1/ h)}{\lambda
  \tau ^{\beta}M(1/ \tau) } \; .
\end{equation}
So we see that for every fixed  real $\kappa \ne 0$  and  positive $s$
\begin{equation}\label{E:eq10}
\widehat{\widetilde{p}}_{h,\tau}(\kappa,s) \,\to \,\frac{s^{\beta
    -1}}{|\kappa|^{\alpha}+ s^{\beta}
}=\widehat{\widetilde{u}}(\kappa,s) \; ,  
\end{equation}
as $h$ and $\tau$ tend to zero under the \textit{scaling relation}
$r(h,\tau)\equiv 1$.
Comparing with (\ref{E:eq6}) we recognize here
$\widehat{\widetilde{u}}(\kappa,s) $  as the combined Fourier-Laplace
transform of the solution to the Cauchy problem
(\ref{E:eq5}). Invoking now the continuity theorems of probability
theory (compare \cite{F1971}) we see that the  time-parameterized  sojourn
probability density converges weakly (or in law) to the solution of
the Cauchy problem (\ref{E:eq5}). We state this result at the
following theorem.
\\
\\
\textbf{Theorem:} Assume the probability laws for the jumps
  $X_{k}$  and the waiting times $T_{k}$
to fulfill the conditions of the Master Lemmata $1$ and $2$,
respectively. Replace the jumps by $hX_{k}$, the waiting times $\tau
T_{k}$. Then for $h$  (and consequently $\tau$) tending 
to zero the solution $p_{h,\tau}(x,t) $
of the rescaled integral equation (\ref{E:eq1}) (the densities there
to be decorated
with indices $h$  and $\tau$)  converges weakly to the solution   of
the Cauchy problem (\ref{E:eq5}), in other words: to the fundamental solution of the
space-time fractional diffusion equation
$\underset{t\;\;
  *}{D}^{\beta}\,u(x,t)=R^{\alpha}\,u(x,t)$.   
%
%%%%%%*****************************************************************
\section{ Examples of random walks}
 Let us first consider the \textit{space-time fractional} diffusion equation 
  more
 closely with regard to special choices of the parameters  $\alpha$
 and $\beta$. In the very particular case $\alpha=2,\; \beta=1$   it
 reduces to the \textit{classical} diffusion equation $\frac{\partial
 u}{\partial t}=\frac{\partial^{2} u}{\partial x^{2}}$.  In the case
 $\alpha=2$, $0<\beta<1$
 we have the \textit{time-fractional} diffusion equation investigated
  in  $1989$ 
  \cite{SW}. In the case $0<\alpha<2$, $\beta=1$  we have the
 \textit{space-fractional} diffusion equation for which the fundamental
 solution is a  symmetric strictly stable probability density
 evolving in time. Whereas the time-fractional case  $\alpha=2$,
  $0<\beta<1$  exhibits \emph{subdiffusive} behaviour, we have
  \emph{superdiffusive} 
  behaviour,  if $0<\alpha<2$. All this can be deduced from the Fourier-Laplace
 representation (\ref{E:eq6})  by observing that the  
 variance  $\langle (x(t))^{2} \rangle = (\sigma(t))^{2}=\int\limits_{-\infty}^{\infty}x^{2}u(x,t)dx$    of the
 position $x(t)$ of a diffusing particle is given as
  $-\frac{\partial^{2}}{\partial\kappa^{2}}\widehat{u}(\kappa,t)|_{\kappa=0}$. Using the
  Mittag-Leffler function  
$E_{\beta}(z)=\sum\limits_{n=0}^{\infty} \frac{z^{n}}{\Gamma(1+n \beta)}$   (see
\cite{GM1997}) we find by Laplace inversion the convergent series
$\widehat{u}(\kappa,t)= E_{\beta}(-|\kappa|^{\alpha} t^{\beta})=1-
\frac{|\kappa|^{\alpha} t^{\beta}}{\Gamma(1+\beta)}+
\frac{|\kappa|^{2\alpha} t^{2\beta}}{\Gamma(1+2\beta)}-+ \cdots $
 from which  for $t>0$ we get $(\sigma(t))^{2}= \frac{2
   t^{\beta}}{\Gamma(1+\beta)} $ if $\alpha=2$, $(\sigma(t))^{2}=
  \infty$ if $0<\alpha<2$.

    Now we present two concrete random walk models, for both assuming
    $0<\alpha<2$, $0<\beta<1$. \emph{The first example} is in
    continuous time and continuous space.   From \cite{ChG} we got the
    idea that it is advantageous to work with functions $W(x)$  and
    $\Phi(t)$  that can elementary be inverted. This is useful to
    produce from $[0,1)$-uniformly distributed pseudo-random numbers
    the jumps and the waiting times in simulations. We take, compare
    \cite{GM2001} and \cite{MRGS}, 
\[
W(x)= \frac{1}{2}+ \frac{1}{2} (-1)^{sign (x)}
\frac{|x|^{\alpha}}{1+|x|^{\alpha}} \; ,\; 
\Phi(t)=1-\frac{ 1}{1+ \Gamma(1-\beta)t^{\beta} } \; .
 \] 
Then we have case (b) of MASTER LEMMA $1$ with  $b=\alpha/2$ and
  $L(x)\equiv 1$,  and case (B) of  MASTER
  LEMMA $2$ with  
  $c=1/|\Gamma(-\beta)|$, $M(t) \equiv 1$.

In the \emph{second example} time and space are equidistantly
discretized. We take $p_{0}=0$,  $p_{k}=b|k|^{-(\alpha+1)}$
for   $0\ne k\, \in \mathbb{Z}$, $c_{n}=cn^{-(\beta+1)}$  for
$n\in \mathbb{N}$, and set $w(x)=\sum\limits_{-\infty}^{\infty}p_{k}
\delta(x-k)$, $\phi(t)=\sum\limits_{k=1}^{\infty}c_{k} \delta(t-k)$
(compare with\cite{GV}). We have case (b) of D-Lemma $1$, case (B) of
D-Lemma $2$ and  
 identify readily (with $\zeta(z)$   denoting Riemann's zeta function)  
$b=\frac{1}{2\zeta(\alpha+1)}$ and $c=\frac{1}{\zeta(\beta+1)}$. The
sequences $p_{1},p_{2},p_{3},\cdots$   and 
$c_{1},c_{2},c_{3},\cdots $     have the nice property of being
completely monotone. The excluded border cases $\alpha=2$  and
$\beta=1$   are singular.

To convey to the reader a feeling for fractional diffusion we present
a few graphical results of approximating random walks, simulated
according to the first example. They show in sequence the case of
Brownian motion (classical diffusion), time-fractional diffusion,
space-fractional diffusion, space-time-fractional diffusion.

  Note that
for $\alpha=2$, we use  the jump density 
 $w(x)=\frac{1}{ \sqrt{2\pi}}\, exp(-\frac{x^{2}}{2})$, for $\beta=1$ the
 waiting time density 
$\phi(t)=exp(-t)$. For  $0<\alpha<2$ and  $0<\beta<1$  we take
the jump distribution  $W(x)$ and
the waiting time distribution $\Phi(t)$, both  as in the first
example.

 Observe some long waiting times in the case $0<\beta<1$ and some long
jumps in the case $0<\alpha<2$.
%
%**************
\begin{figure}[ht!]
\begin{minipage}{4cm}
%\centering
\epsfysize=50mm
\epsfbox[50 360 500 580]{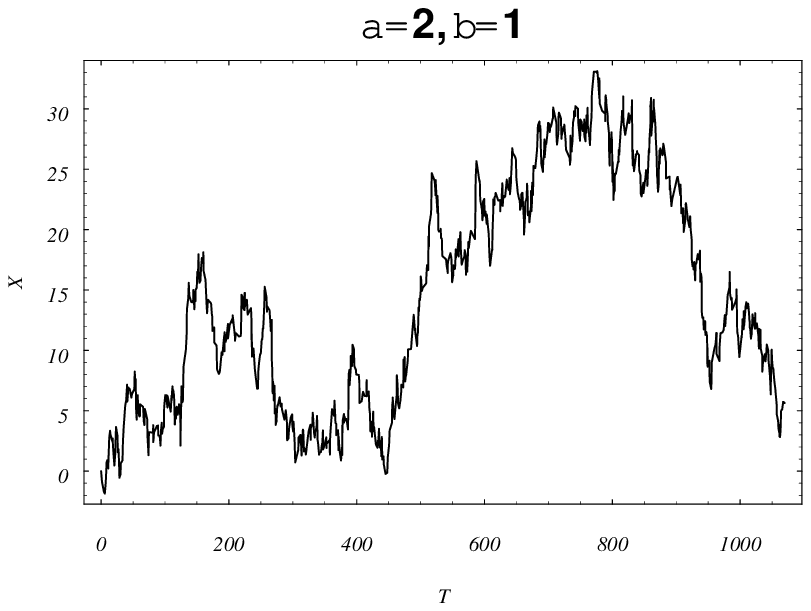}
\caption{normal diffusion}
\end{minipage}
\hspace*{2cm}
\begin{minipage}{4cm}
%\centering
\epsfysize=50mm
\epsfbox[50 360 500 580]{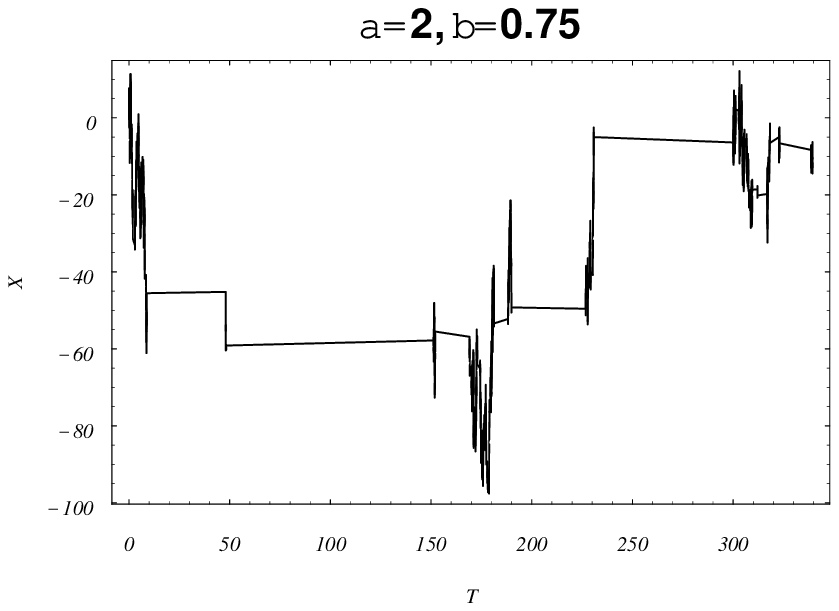}
\caption{time-fractional}
\end{minipage}
\end{figure}
\begin{figure}[ht!]
\begin{minipage}{4cm}
%\centering
\epsfysize=50mm
\epsfbox[50 360 500 580]{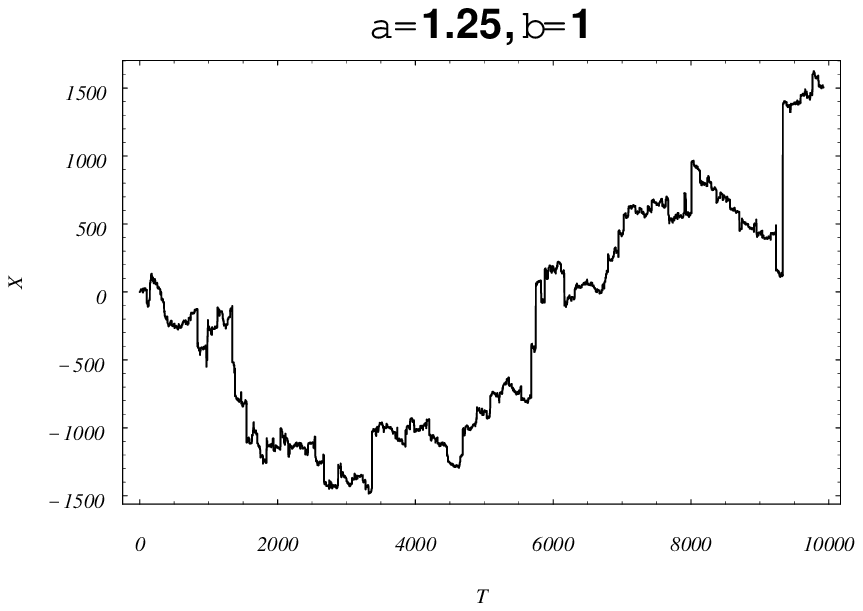}
\caption{space-fractional}
\end{minipage}
\hspace*{2cm}
\begin{minipage}{4cm}
%\centering
\epsfysize=50mm
\epsfbox[50 360 500 580]{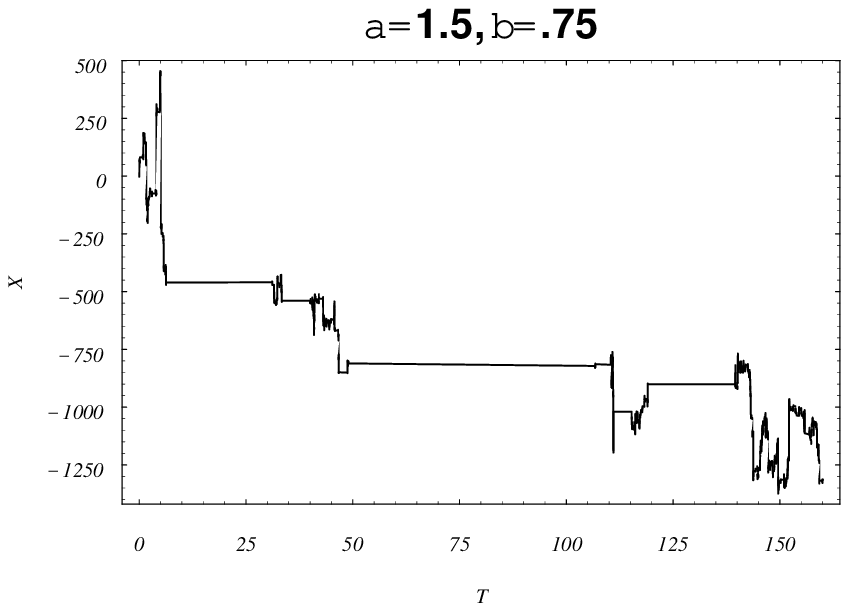}
\caption{space-time fractional}
\end{minipage}
\end{figure}
\pagebreak
%*********************************************
\section{Comments, suggestions and conclusions}
     The theory of compound renewal processes, also called renewal
     processes with reward, in physics and other natural sciences
     called  \emph{continuous time random walks} (though space and time
     need not be continuous) began to flourish in the middle of the
     sixties of the past century, let us quote \cite{MW} and
     \cite{C}. We cannot give here a comprehensive survey of relevant
     literature, so we ask all not mentioned contributors to forgive
     us this surely biased account. For larger lists of references and
     more competent appreciation of achievements and applications we
     recommend \cite{Bar} and \cite{MK}. As an early pioneer
     Balakrishnan  \cite{Bal} deserves to be put into  light. He
     has, in $1985$, found the time-fractional diffusion equation
     ($\alpha=2$ , $0<\beta<1$) as the properly scaled diffusion limit
     for some random walks with power law waiting time. At that time,
     four years before in \cite{SW} the basic analytic theory was 
     developed, the name \emph{fractional diffusion} was not yet
     common, and so was not used in \cite{Bal},  hence Balakrishnan
     did not find the resonance he would have deserved. A decisive
     step forward occurred in \cite{HA} in $1995$.
 Hilfer and Anton, roughly speaking, showed among other things that by
     taking the Mittag-Leffler waiting time density
     $-\frac{d}{dt}E_{\beta}(-t^{\beta})$    the basic  equation of
     continuous time random walk can be transformed to a 
     time-fractional evolution equation for the sojourn probability
     density. Thus they have essentially  found the time-fractional
     generalization 
     of the Kolmogorov-Feller evolution equation for the compound
     Poisson process which e.g. is treated in \cite{F1971}.  However,
     already in \cite{Bal} appears the waiting time density whose
     Laplace transform is $(1+s^{\beta})^{-1}$  as playing a
     distinct role, but was not recognized as a function of
     Mittag-Leffler type (such functions too long having been
     insufficiently known). Gorenflo and Mainardi and co-authors have,
     beginning in $1998$, published several papers on various types of
     approximating random walks for space-fractional and space-time
     fractional diffusion processes of which we quote \cite{GM2001},
     \cite{GM2002}, \cite{GM2003}  and \cite{GV}, furthermore some
     papers (stressing the relevance of the Mittag-Leffler waiting
     time) motivated by applications to finance: \cite{SGM},
     \cite{MRGS}, \cite{GMSR}. In \cite{GMSR} the space-time fractional
     diffusion equation is obtained as a diffusion limit of the
     time-fractionalized 
     Kolmogorov-Feller equation
\begin{equation}\label{E:eq17}
\underset{t\;\;
  *}{D}^{\beta}\,p(x,t)= -p(x,t)+\int\limits_{-\infty}^{\infty}
w(x-x') p(x',t) dx' \; .
\end{equation}
 The publications \cite{MLP} and
     \cite{MPG} are devoted to analytic treatment via integral
     representations of the evolving probability densities that solve
     the (spatially one-dimensional) space-time fractional diffusion
     equation.  An important concept in fractional diffusion processes
     is the concept of subordination (see, e.g. \cite{MS} and
     \cite{MBSB}). By our way of relating the scaling parameters in
     the passage to the limit we circumvent this concept.   
     Let us mention again \cite{GM2002}. There we have based our
     scaled transition to the limit on two lemmata for the asymptotics
     of the transforms of the densities whereas here we work with the
     MASTER LEMMATA for the distribution functions, motivated by
     \cite{GK}. This, of course, is more general and allows discrete
     and continuous probabilities and mixtures of them. However, in
     not so general situations the other lemmata may be simpler to
     apply (as we have done in \cite{GV} for the fully discrete case
     with regular grids). In \cite{US}, in contrast to our treatment in
     \cite{GM2002} and here, the scaling is not done via the individual
     steps in space and time but directly in the distributions of
     waiting times and jumps. However, this is equivalent to our
     way. Let us in this context say a few words to the essential
     statement  of \cite{H}. There Hilfer shows that a power law for
     the waiting time is not sufficient for getting in the limit a
     fractional diffusion process with the fractional time derivative
     having the same order as the power law.  This seemingly negative
     result, however, does not hit the theory expanded here in our
     paper. In Hilfer's counter-example the passage to the diffusion
     limit is not well-scaled in our sense; in fact, in it are hidden
     two different scalings. Thus \cite{H} may inspire to investigate
     systematically continuous time random walks that can be scaled in
     more than one way.

Let us, as a final statement, say that the case of non-symmetric jump
distributions (bypassed in our paper) can analogously be studied, and
let us also hint to \cite{Re}.
\\
\\
%+++++++++++++++++++++++++++++++++++++++++++++++++++++++++++++++++
\textbf{Acknowledgments} 
This work has partially been carried out in the frame of the  
 \textit{INTAS} project 00-0847. The second named author is 
grateful for  the grant provided by the government of \textit{Arab Republic of
Egypt}.  We are grateful to  F. Mainardi and E. Scalas for 
fruitful discussions on the subject. The first named author  thanks
R. Hilfer for a preprint of \cite{H} and for inspiring discussions. 
%++++++++++++++++++++++++++++++++++++++++++++++++++++++++++++++++++

%
\end{document}